\documentclass[amsthm]{elsart}
\usepackage{ifpdf}
\usepackage{amssymb,lineno}
\usepackage{float}
 \usepackage[dvipdfm]{graphicx}
\usepackage{tikz}
\usepackage{array}
\usepackage{color}
\usepackage{mathrsfs}
\usetikzlibrary{arrows}

\ifpdf
\usepackage[%
  pdftitle={Instructions for use of the document class
    elsart},%
  pdfauthor={Simon Pepping},%
  pdfsubject={The preprint document class elsart},%
  pdfkeywords={instructions for use, elsart, document class},%
  pdfstartview=FitH,%
  bookmarks=true,%
  bookmarksopen=true,%
  breaklinks=true,%
  colorlinks=true,%
  linkcolor=blue,anchorcolor=blue,%
  citecolor=blue,filecolor=blue,%
  menucolor=blue,pagecolor=blue,%
  urlcolor=blue]{hyperref}
\else
\usepackage[%
  breaklinks=true,%
  colorlinks=true,%
  linkcolor=blue,anchorcolor=blue,%
  citecolor=blue,filecolor=blue,%
  menucolor=blue,pagecolor=blue,%
  urlcolor=blue]{hyperref}
\fi

\journal{~~}

\makeatletter
\def\elsartstyle{%
    \def\normalsize{\@setfontsize\normalsize\@xiipt{14.5}}
    \def\small{\@setfontsize\small\@xipt{13.6}}
    \let\footnotesize=\small
    \def\large{\@setfontsize\large\@xivpt{18}}
    \def\Large{\@setfontsize\Large\@xviipt{22}}
    \skip\@mpfootins = 18\p@ \@plus 2\p@
    \normalsize
}
\@ifundefined{square}{}{\let\Box\square}
\makeatother

\begin{document}

\begin{frontmatter}
\title{The characterization of perfect Roman  domination stable trees}

\author[LanZhou]{Zepeng Li}\ead{lizp@lzu.edu.cn},
\author[GZHU]{Zehui Shao}\ead{zshao@gzhu.edu.cn},
\author[GZHU]{Yongsheng Rao}\ead{rysheng@gzhu.edu.cn},
\author[CDU]{Pu Wu}\ead{puwu1997@126.com},
\author[Adelphi]{Shaohui Wang}\ead{shaohuiwang@yahoo.com}

\address[LanZhou]{
    School of Information Science and Engineering, Lanzhou University, Lanzhou 730000, China}
\address[GZHU]{
   Institute of Computing Science and Technology, Guangzhou University, Guangzhou 510006, China}
\address[CDU]{ School of Information Science and Engineering,
 Chengdu University,
  Chengdu 610106,  China}
\address[Adelphi]{Department of Mathematics,
Savannah State University, Savannah, GA 31419, USA}
\corauth[cor]{Corresponding author: Shaohui Wang.}

\begin{abstract}
A \emph{perfect Roman dominating function} (PRDF) on a graph $G = (V, E)$ is a function $f :
V \rightarrow \{0, 1, 2\}$ satisfying the condition that every vertex $u$ for which $f(u) = 0$
is adjacent to exactly one vertex $v$ for which $f(v) = 2$.
The weight of a PRDF is the value $w(f) = \sum_{u \in V}f(u)$. The minimum weight of
a PRDF on a graph $G$ is called the \emph{perfect Roman domination
number $\gamma_R^p(G)$} of $G$.
A graph $G$ is perfect Roman domination domination stable if the perfect Roman domination number of $G$ remains unchanged
under the removal of any vertex.
In this paper, we characterize  all trees that are perfect Roman domination stable.

\end{abstract}

\end{frontmatter}

\newtheorem{thm1}{Theorem}
\newtheorem{clm1}{Claim}
\newtheorem{rem1}{Remark}
\newtheorem{lem1}{Lemma}
\newtheorem{cor1}{Corollary}
\newtheorem{def1}{Definition}
\newtheorem{observ1}{Observation}
\newtheorem{ex1}{Example}
\newtheorem{prop1}{Proposition}

\section{Introduction}
In this paper, we shall only consider graphs without multiple edges or loops.
Let $G$ be a graph, $S \subseteq V(G)$, $v \in V(G)$, the \emph{neighborhood} of $v$ in $S$ is denoted by $N_S(v)$.
That is to say $N_S(v) =\{u| uv \in E(G), u \in S\}$.
The closed neighborhood $N_S[v]$ of $v$ in $S$ is defined as $N_S[v]=\{v\} \cup N_S(v)$.
If $S=V(G)$, then $N_S(v)$ and $N_S[v]$ are denoted by $N(v)$ and $N[v]$, respectively.
Let $S \subseteq V(G)$, we write $N_G(S)= \cup_{x \in S}N_G(x)$.
The degree of $v$ is $d(v)=|N(v)|$.
We will omit the subscript $G$, that is to say, $N_G(T)$ is denoted by $N(T)$.
For a tree $T$ and a vertex $v\in V(G)$, we denote by $L(v)$ the set of all leaves of $v$.
 A tree $T$ is a double star if it contains exactly two vertices that are not leaves.
A double star with respectively $p$ and $q$ leaves attached at each support vertex is denoted by $DS_{p,q}$.

A \emph{dominating set} of $G$ is a subset $D$ of $V$ such that every vertex in $V- D$ is adjacent to at least one vertex in $D$. The \emph{domination number} $\gamma(G)$ is the minimum cardinality of a dominating set of $G$. A \emph{Roman dominating function} (RDF) of $G$ is a function $f: V(G)\rightarrow \{0, 1, 2\}$ such that every vertex $v\in V$ with $f(v)=0$ is adjacent to at least one vertex $u$ with $f(u)=2$. The weight of an RDF $f$ is the value $f(V)=\sum_{v\in V} f(v)$. The Roman
domination number $\gamma_R(G)$ is the minimum weight of an RDF of $G$. The problems on domination and Roman domination of graphs have been investigated widely, for example, see the list of references, \cite{Haynes1998} and \cite{Chambers2009,Cockayne2004,Favaron2009,WW6,WW7,WW8}, respectively.

The affections of vertex removal on domination number and Roman domination number in a graph
have been studied in \cite{Bauer,Rad1} and \cite{Hansberg,Rad2}, respectively. Jafari Rad and Volkmann \cite{Rad2} introduced the concept
of Roman domination stable graphs and these graphs had been further studied in \cite{Hajian,Samodivkin,WW1,WW2,WW3,WW4,WW5}.

Furthermore, Henning, Klostermeyer and MacGillivray\cite{Henning} introduce a perfect version of Roman
domination.
A \emph{perfect Roman dominating function} (PRDF) on a graph $G = (V, E)$ is a function $f :
V \rightarrow \{0, 1, 2\}$ satisfying the condition that every vertex $u$ for which $f(u) = 0$
is adjacent to exactly one vertex $v$ for which $f(v) = 2$.
The weight of a PRDF is the value $w(f) = \sum_{u \in V}f(u)$. The minimum weight of
a PRDF on a graph $G$ is called the \emph{perfect Roman domination
number $\gamma_R^p(G)$} of $G$.
A PRDF $f$ is call a $\gamma_R^p$-function of $G$ if $w(f)=\gamma_R^p(G)$.
A graph $G$ is \emph{perfect Roman domination stable} if the perfect Roman domination number of $G$ remains unchanged
under removal of any vertex.

Recently, many research are working in this topic. For instance,
Rad et al.\cite{Rad2} studied the changing and unchanging the Roman domination number of a graph. Henning et al. \cite{Henning} explored some trees about perfect Roman domination. Favaron et al. \cite{Favaron2009} found some Romain domination number of a graph. Chambers et al. \cite{Chambers2009} deduced some extremal results on Roman domination.
Motivated by the above results,  we continue to study Roman domination and characterize all perfect Roman domination stable trees.

\section{Perfect Roman domination stable trees}
In this section, we will give some lemmas and properties.
\begin{observ1}\label{PRO1}
Let $T$ be a perfect Roman domination stable tree and $f$ be a $\gamma_R^p$-function of $T$. Then \\
(i)  $f(v)\neq 1$ for any $v \in V(T)$.\\
(ii) $f(v) \neq 2$ for any leaf $v$ in $T$.\\
(iii) If there exists a vertex $x_4 \in V(T)$ is adjacent to a star with vertex set $\{x_1,x_2,x_3,y_1\}$ and center $x_2$ and $x_4x_3 \in E(T)$,
then $f(x_2)=2$ and  $f(x_1)=f(y_1)=f(x_3)=f(x_4)=0$.
\end{observ1}
\begin{proof}
(i) Suppose there exists a vertex $v\in V(T)$ for which $f(v)=1$.
 Let $T'=T-\{v\}$. Since $T$ is a perfect Roman domination stable tree, we have $\gamma_R^p(T)=\gamma_R^p(T')$.
Then we have $f|_{T'}$ is a PRDF on $T'$.
Thus
$\gamma_R^p(T')\leq w(f)-1= \gamma_R^p(T)-1$
a contradiction.\\
(ii)
Otherwise, let $u \in N(v)$ and $f(u)=0$.
Now we can obtain a PRDF with the same weight by assigning 1 to $u$ and $v$, a contradiction with (i).\\
(iii) By the results of (i) and (ii), we have $f(x_2)=2$, $f(x_1)=f(y_1)=0$ and $f(x_3) \in \{0,2\}$.
If $f(x_3)=2$, then we can obtain a new $\gamma_R^p$-function of $T$ by changing $f(x_2)$ to 0 and $f(s)$ to 1 for any $s\in \{x_1,y_1\}$,
contradicting with (i).
Thus, $f(x_3)=f(x_4)=0$.
\end{proof}

In this section we give a constructive characterization of all perfect Roman domination stable trees under vertex removal.
For a tree $T$, let
\begin{eqnarray*}
W(T) & = & \{u \in V(T)\mid \textrm{ $f(u)=0$ for any $\gamma_R^p$-funnction $f$ of $T$}\}.
\end{eqnarray*}

In order to presenting our constructive characterization, we define a family of trees as follows.
Let $\mathcal{T}$ be the family of trees $T$ that can be obtained from a sequence $T_1$,
$T_2$, $\ldots$, $T_k$ of trees for some $k \geq 1$, where $T_1$ is $P_3$ and $T=T_k$.
If $k \geq 2$, $T_{i+1}$ can be obtained from $T_i$ by    the following operation.

\begin{description}
  \item [Operation ${\mathcal O}_1$: ] If $u \in W(T_i)$, then ${\mathcal
O}_1$ adds a path $v_3v_2v_1$ and an edge $uv_3$  to obtain $T_{i+1}$.
\end{description}

\begin{lem1}\label{AddP3}
Let $G$ be a graph and $u \in V(G)$.
If $G'$ is a graph obtained by adding a  path  $v_3v_2v_1$ and an edge $uv_3$ from $G$,
then $\gamma_R^p(G')= \gamma_R^p(G) +2$.
\end{lem1}

\begin{proof}
Let $f$ be a $\gamma_R^p$-function on $G$.
If $f(u)\neq 2$, we define $f' :V(G') \rightarrow \{0, 1, 2\}$ by $f'(v_1)=0$, $f'(v_2)=2$, $f'(v_3)=0$ and $f'(v)=f(v)$ if $v\in V(G)$.
If $f(u)=2$, we
define $f' :V(G') \rightarrow \{0, 1, 2\}$ by $f'(v_1)=1$, $f'(v_2)=1$, $f'(v_3)=0$ and $f'(v)=f(v)$ if $v\in V(G)$.
Then in each case $f'$ is a PRDF function on $G'$ and $w(f')=w(f)+2$. Thus we have
$\gamma_R^p(G')\leq w(f')=w(f)+2 = \gamma_R^p(G)+2$.

Conversely, let $f$ be a $\gamma_R^p$-function on $G'$.  We consider the following cases.\\
\textbf{Case 1:} $f(u) \leq 1$. Clearly, we have $f(v_2)=2$ and $f(v_1)=f(v_3)=0$. Then $f|_G$ is a PRDF function on $G$.
Thus $\gamma_R^p(G')= w(f)= w(f|_G)+2 \geq \gamma_R^p(G)+2$, as desired.\\
\textbf{Case 2:}  $f(u)=2$. In this case we have $(f(v_3),f(v_2),f(v_1)) \in \{(0,1,1),(0,0,2)\}$.
 Now $f|_G$ is a PRDF of $G$ and so $\gamma_R^p(G) \leq w(f|_G)=w(f)-2=\gamma_R^p(G')-2$, as desired.
\end{proof}

\begin{lem1}\label{OneV}
Let $T$ be a perfect Roman domination stable tree and $u \in W(T)$.
If $T'$ is a tree obtained by adding a single vertex $v$ and an edge $uv$ from $T$,
then $\gamma_R^p(T')= \gamma_R^p(T) +1$.
\end{lem1}
\begin{proof}
Let $f$ be a $\gamma_R^p$-function on $T$. Define $f' :V(T') \rightarrow \{0, 1, 2\}$ by $f'(v)=1$ and $f'(x)=f(x)$ if $x\in V(T)$. Then $f'$ is a PRDF function on $T'$ and $w(f')=w(f)+1$. Thus we have
$\gamma_R^p(T')\leq w(f')=w(f)+1 \leq \gamma_R^p(T)+1$.

Conversely, let $f$ be a $\gamma_R^p$-function on $T'$, by Observation \ref{PRO1}  we have $f(v)\neq 1$ for any $v \in V(T)$.
Now we consider the following cases. \\
\textbf{Case 1:} $f(v)=2$. In this case we have $f(u) = 0$ and we consider a $\gamma_R^p$-function $f'$ on  $T'$ with $f'(v)=f'(u)=1$ and $f'(x)=f(x)$ otherwise. Then  $f'|_T$ is a PRDF of $T$ and we have $\gamma_R^p(T)\leq \gamma_R^p(T')-1$, as desired.\\
\textbf{Case 2:} $f(v)=0$. In this case we have $f(u) =2$.
Since $u \in W(T)$, and  $f|_T$ is a PRDF of $T$.
Then we have $\gamma_R^p(T)< w(f|_T)=\gamma_R^p(T')$.
Thus $\gamma_R^p(T)\leq \gamma_R^p(T')-1$, as desired.
\end{proof}

\begin{lem1}\label{AddP2}
Let $T$ be a perfect Roman domination stable tree and $u \in W(T)$.
If $T'$ is a tree obtained by adding a  path  $v_2v_1$ and an edge $uv_2$ from $T$,
then $\gamma_R^p(T')= \gamma_R^p(T) +2$.
\end{lem1}
\begin{proof}
Let $f$ be a $\gamma_R^p$-function on $T$. Define $f' :V(T') \rightarrow \{0, 1, 2\}$ by $f'(v_2)=0$, $f'(v_1)=2$ and $f'(v)=f(v)$ if $v\in V(T)$. Then $f'$ is a PRDF function on $T'$ and $w(f')=w(f)+2$. Thus we have
$\gamma_R^p(T')\leq w(f')=w(f)+2 = \gamma_R^p(T)+2$.

Conversely, let $f$ be a $\gamma_R^p$-function on $T'$, by Observation \ref{PRO1}  we have  $f(v)=0$ for any $v \in V(T)$.
Now we consider the following cases. \\
\textbf{Case 1:} $f(v_1)=2$. In this case we have $f(v_2) = 0$.   Then $f|_T$ is a PRDF of $T$ and we have $\gamma_R^p(T)\leq \gamma_R^p(T')-2$, as desired.\\
\textbf{Case 2:} $f(v_1)=0$. In this case we have $f(v_2) =2$ and $f(u)=0$.
Since $T$ is stable,
we have $\gamma_R^p(T-u)=\gamma_R^p(T)$.
Since $f(u)=0$, we have $f|_{T-u}$ is a PRDF of $T-u$, we have
$\gamma_R^p(T)=\gamma_R^p(T-u) \leq w(f|_{T-u})=w(f)-2=\gamma_R^p(T')-2$, as desired.
\end{proof}

\begin{lem1}\label{first}
{\em If $T_i$ is a perfect Roman domination stable tree and $T_{i+1}$ is a tree obtained from
$T_i$ by Operation ${\mathcal O}_1$, then
$T_{i+1}$ is a perfect Roman domination stable tree.}
\end{lem1}
\begin{proof}
By Lemma \ref{AddP3}, we have $\gamma_R^p(T_{i+1})=\gamma_R^p(T_{i})+2$.
Let $v \in V(T_{i+1})$ be an arbitrary vertex and  $T'=T_{i+1}-v$. \\
If $v\in V(T_{i})\setminus\{u\}$,
since $T_i$ is  stable, we have $\gamma_R^p(T_i-v)=\gamma_R^p(T_i)$.
Then by Lemma \ref{AddP3},  we have $$\gamma_R^p(T_{i+1}-v)=\gamma_R^p(T_i-v)+2=\gamma_R^p(T_i)+2=\gamma_R^p(T_{i+1}).$$\\
If $v=u$, since $T_i$ is stable,
obviously, $\gamma_R^p(T')=\gamma_R^p(T_i-v)+2=\gamma_R^p(T_i)+2=\gamma_R^p(T_{i+1})$.\\
If $v=v_1$, By Lemma \ref{AddP2}, we have $\gamma_R^p(T')=\gamma_R^p(T_i)+2=\gamma_R^p(T_{i+1})$.\\
If $v=v_2$, By Lemma \ref{OneV}, we have $\gamma_R^p(T')=\gamma_R^p(T_i+{uv_3})+\gamma_R^p(K_1)=\gamma_R^p(T_i)+1+1=\gamma_R^p(T_{i+1})$.\\
If $v=v_3$, then we have $\gamma_R^p(T')=\gamma_R^p(T_i)+\gamma_R^p(P_2)=\gamma_R^p(T_i)+2=\gamma_R^p(T_{i+1})$.
Thus $T_{i+1}$ is a perfect Roman  domination stable and the proof is complete.
\end{proof}

\begin{lem1}\label{Lemma5}
Let $T$ be a perfect Roman domination stable tree  of order $n\ge 3$ with $diam(T) \geq 4$ and $T$ contain no pendent $P_3$, and $P=x_1x_2\cdots x_k$ be a longest path of $T$.
Then $d(x_2)=2$.
\end{lem1}
\begin{proof}
Let $f$ be a $\gamma_R^p$-function of $T$.
By Observation \ref{PRO1}, we have $f(v) \in \{0,2\}$ for any $v \in V(T)$.

\begin{figure}[H]
\begin{center}
\includegraphics[scale=0.20]{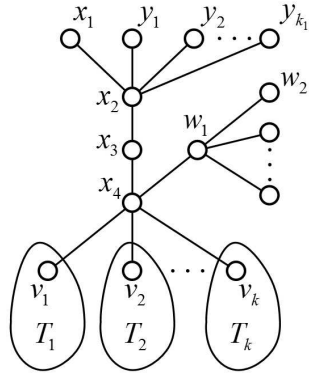}
\caption{The tree in the proof of Lemma \ref{Lemma5}}\label{LV2}
\end{center}
\end{figure}

%
First we have $|L(x_2)| \leq 2$. Otherwise, we assume $L(x_2)=\{x_1, y_1, \cdots, y_{k_1}\}$ where $k_1 \geq 2$.
It is clear that $f(x_2)=2$, $f(x_1)=0$ and $f(y_i)=0$ for any $i=1,2,\cdots,k_1$.
Let $g$ be a $\gamma_R^p$-function of $T-x_2$ and so $g(x_1)=g(y_i)=1$ for any $i=1,2,\cdots,k_1$.\\
If $g(x_3)=2$,  we define a function $g'$ of $T$ as follows.
$g'(x_2)=2$, $g'(x_1)=g'(y_i)=0$ for any $i=1,2,\cdots,k_1$ and $g'(v)=g(v)$ for other vertex $v \in T$ and we have
$\gamma_R^p(T)\leq w(g')=w(g)-(k_1+1)+2<w(g)=\gamma_R^p(T-x_2)$, a contradiction.\\
If $g(x_3)=0$,  we define a function $g'$ of $T$ as follows.
$g'(x_2)=2$, $g'(x_3)=1$, $g'(x_1)=g'(y_i)=0$ for any $i=1,2,\cdots,k_1$ and $g'(v)=g(v)$ for other vertex $v \in T$. If $k_1\geq 3$, then we have $\gamma_R^p(T)\leq w(g')=w(g)-(k_1+1)+3<w(g)=\gamma_R^p(T-x_2)$, a contradiction. If $k_1=2$, then we have $w(g')=w(g)$ but $x_3$ is assigned with 1, a contradiction with Observation \ref{PRO1}.

Then we have $|L(x_2)| = 1$.
Suppose to the contrary, $|L(x_2)| = 2$, then we have
\begin{clm1}
$d(x_3)=2$.
\end{clm1}
\begin{proof}
First, we have $x_3$ has no leaf neighbor.
Otherwise, let $w_1 \in L(x_3)$.
Since $f(x_3)=0$ from Observation \ref{PRO1} (iii), we have $f(w_1)=1$, contradicting Observation \ref{PRO1}.

Now we have $x_3$ is not adjacent to a pendent star.
Otherwise, assume $x_3$ is adjacent to a pendent star centered with $w_1$ with $x_3w_1 \in E(T)$ and $L(W_1)$.
Since $f(x_3)=0$ and $f(x_2)=2$ from Observation \ref{PRO1} (iii), we have $f(w_1)=0$, and $f(x)=1$ for any vertex $v\in L(w_1)$, contradicting with Observation \ref{PRO1}.
\end{proof}

Now we have $x_4$ has no leaf neighbor.
Otherwise, let $w_1 \in L(x_4)$.
Since $f(x_3)=0$ and $f(x_2)=2$ , we have $f(x_4)=0$ and thus $f(w_1)=1$, contradicting Observation \ref{PRO1}.
Now  it follows from Observation \ref{PRO1}  that $f(x_1)=f(x_3)=f(y_1)=f(x_4)=0$, $f(x_2)=2$ and we have
\begin{clm1}\label{Claim1}
There exists no vertex $w_1 \in N(x_4)$ for which $|L(w_1)|=d(w_1)-1 \geq 1$.
\end{clm1}
\begin{proof}
Otherwise, we assume $w_1 \in N(x_4)$ for which $|L(w_1)|=d(w_1)-1 \geq 1$.
Analogous to the proof of $|L(x_2)| \leq 2$, we have $|L(w_1)| \leq 2$.
Let $T'=T-\{x_1,y_1,x_2,x_3\}-N[w_1]$ with $T'=T_1 \cup T_2 \cup \cdots \cup T_k$ for some $k$ and
$v_i = N(x_4) \cap T_i$ for $i=1,2,\cdots,k$ (see Fig. \ref{LV2}).
By Observation \ref{PRO1}, we have $f(x_2)=2$, $f(x_3)=f(x_4)=0$, $f(w_1)=2$ and $f(v_i)=0$ for any $i$.\\

Then we have
\begin{equation}\label{EQAdd4}
 \gamma_R^p(T)=w(f)= \sum_{i=1}^k w(f|_{T_i}) + 4.
\end{equation}

Since $f|_{T_i}$ is a PRDF of $T_i$,   we have $\gamma^p_R(T_i) \leq w(f|_{T_i})$ for any $i$.

Now we will show that
\begin{equation}\label{EQ1}
\textrm{$\gamma^p_R(T_i) = w(f|_{T_i})$ for any $i$.}
\end{equation}
Otherwise, we assume that $g$ is a $\gamma_R^p$-function of $T'$.
Now we suppose $\gamma^p_R(T_i) < w(f|_{T_i})$, then we have $w(g|_{T_i}) < w(f|_{T_i})$ for some $i$.
If $g(v_i)\leq 1$, then under $f$ we assign $g(t)$ instead of $f(t)$ to $t$ for any $t \in V(T_i)$ and obtain a PRDF with fewer weight, a contradiction.
If $g(v_i)=2$,  then under $f$ we assign  $g(t)$ instead of $f(t)$ to $t$ with for any $t \in V(T_i)$, $g(z)=1$ for $z \in N[w_1]-x_4$ and obtain a PRDF with weight at most $w(f)$. But $w_1$ is assigned with 1,  contradicting Observation \ref{PRO1}.

Since $f|_{T_i-v_i}$ is a PRDF of $T_i-v_i$,   we have $\gamma^p_R(T_i-v_i) \leq w(f|_{T_i-v_i})=w(f|_{T_i})$ for any $i$.

Now we will show that
\begin{equation}\label{EQ2}
\textrm{$\gamma^p_R(T_i-v_i) = w(f|_{T_i-v_i})=w(f|_{T_i})$ for any $i$.}
\end{equation}
Otherwise, we
assume that $g$ is a $\gamma_R^p$-function of $T'$.
and suppose $\gamma^p_R(T_i-v_i) < w(f|_{T_i})$. Then we have $w(g|_{T_i-v_i}) < w(f|_{T_i-v_i})$ for some $i$.
Then under $f$ we assign 1 to $v_i$  and $g(t)$ instead of $f(t)$ to $t$ for any $t \in V(T_i-v_i)$ and obtain a PRDF of $T$ with at most weight $w(f)$.
But $v_i$ is assigned with 1, contradicting Observation \ref{PRO1}.

Now let $g$ be a $\gamma_R^p$-function of $T-x_2$.
Now we have $\sum\limits_{z \in N[w_1]} g(z) \geq 3$.
Since $T$ is stable, we have $w(f)=w(g)$.
Therefore, we have
\begin{equation}\label{4}
\gamma_R^p(T)=w(g)=\sum_{i=1}^k w(g|_{T_i}) + 2+ \sum_{z \in N[w_1]} g(z) \geq \sum_{i=1}^k w(g|_{T_i})+5,
\end{equation}
By Eqs.(\ref{4}) and (\ref{EQAdd4}), we have
\begin{equation}\label{EQAdd5}
 \sum_{i=1}^k w(g|_{T_i})\leq \gamma_R^p(T)-5< \sum_{i=1}^k w(f|_{T_i}).
\end{equation}
Hence by Eq.(\ref{EQ1}), there must exist some $j$ ($1\leq j \leq k$) satisfying that
\begin{equation}\label{EQAdd6}
w(g|_{T_j})<w(f|_{T_j})=\gamma_R^p(T_j).
\end{equation}
If $g(v_j)\geq 1$, then $g|_{T_j}$ is a PRDF on $T_j$, a contradiction with Eq.(\ref{EQAdd6}).
If $g(v_j)=0$, then $g|_{T_j-v_j}$ is a PRDF on $T_j-v_j$. By Eq.(\ref{EQ2}), $w(g|_{T_j-v_j})=w(g|_{T_j})<w(f|_{T_j})=\gamma_R^p(T_j-v_j)$, a contradiction.

This completes the proof.
\end{proof}

\begin{clm1}\label{Claim2}
$x_4$ is not adjacent to a vertex $w_1$ with $|L(w_1)|\geq 2$.
\end{clm1}
\begin{proof}
Suppose to the contrary, there exists  a vertex $w_1\in N(x_4)$  with $|L(w_1)|\geq 2$.
By Claim \ref{Claim1}, there exists a vertex $w_2\in N(w_1)\setminus \{x_4\}$ such that $|L(w_2)|\leq 2$.
We deduce from Observation \ref{PRO1} that $f(x_4)=f(y)=0$ for any $y\in L(w_1)\cup L(w_2)$, and $f(w_1)=f(w_2)=2$.

Then we consider a PRDF $f'$ of $T$ with $f'(w_2)=0$ and $f'(x)=1$ for any $x\in L(w_2)$.
It is obvious that $w(f')\leq w(f)$ and $f'$ is a $\gamma_R^p$-function of $T$.
But there exists a leaf vertex which is assigned with $1$ under $f'$, a contradiction with Observation \ref{PRO1}.

\end{proof}

\begin{clm1}
 $d(x_4)=2$.
\end{clm1}
\begin{proof}
Suppose to the contrary $d(x_4) \geq 3$, then by Claims \ref{Claim1} and \ref{Claim2}, it is sufficient to consider the following four cases.

\begin{figure}[H]
\begin{center}
\includegraphics[scale=0.7]{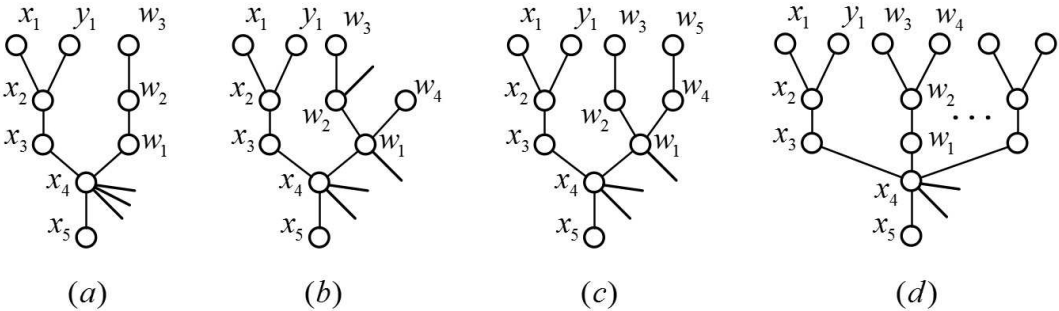}
\caption{Four cases if $d(x_4) \geq 3$ in the proof of  Lemma \ref{Lemma5}}
\label{FourCases}
\end{center}
\end{figure}

\textbf{Case A:} $x_4$ is adjacent to a pendent $P_3$ (see Fig. \ref{FourCases}a), contradicting with the condition of $T$.\\
\textbf{Case B:} $x_4$ is adjacent to a $P_4=w_3w_2w_1w_4$ with $x_4w_1 \in E(T)$ (see Fig. \ref{FourCases}b).\\
  Easily, we have $f(w_1)=f(w_3)=0$ and $f(w_2)=2$. Hence $f(w_4)=1$, a contradiction.

\textbf{Case C:} $x_4$ is adjacent to a $P_5=w_3w_2w_1w_4w_5$ with $x_4w_1 \in E(T)$ (see Fig. \ref{FourCases}c).\\
By Observation \ref{PRO1}, $f(w_3)=f(w_5)=0$ and $f(w_2)=f(w_4)=2$.
If $f(w_1)=0$, a contradiction (a vertex in $V_0$ is adjacent to exactly one vertex in $V_2$). If $f(w_1)=2$, then under $f$ we assign
   2,0,1,0,1 to $w_1,w_2,w_3,w_4$ and $w_5$, respectively,
and obtain a PRDF of $T$ with fewer weight $w(f)$, a contradiction.

\textbf{Case D:} $x_4$ is adjacent to a star with vertex set $\{w_1,w_2,w_3,w_4\}$ centered at $w_2$ with $x_4w_1 \in E(T)$ (see Fig. \ref{FourCases}d).\\
We can obtain  $f(w_1)=f(w_3)=f(w_4)=0$ and $f(w_2)=2$.

Let $T'=T-\{w_2,w_3,w_4\}$. Let $g$ be a $\gamma_R^p$-function on $T'$.
Since $T$ is stable,  We have $w(f)=\gamma_R^p(T-w_2)=w(g)+2$.

If $g(w_1)=2$, under $g$ we assign  0  to $w_2$,  1 to $w_3,w_4$ and obtain a PRDF of $T$ with same weight $w(f)$. But $w_3$ and $w_4$ are assigned with 1, contradicting Observation \ref{PRO1}.

If $g(w_1)=1$, under $g$ we assign  2  to $w_2$,  0 to $w_3,w_4$ and obtain a PRDF of $T$ with same weight $w(f)$. But $w_1$ is assigned with 1, contradicting Observation \ref{PRO1}.

If $g(w_1)=0$, then $g(x_4)=2$. Hence we have $g(x_1)=g(x_2)=g(y_1)=1$ and $g(x_3)=0$.
Thus $g(x_5)=0$. Otherwise, under $g$ we assign 2 to $x_2,w_2$, 0 to $x_1,y_1,x_3,w_1,w_3,w_4$, 1 to $x_4$ and obtain a PRDF of $T$ with fewer weight $w(f)$, a contradiction.

Note that $g(x_4)=2$ and $g(x_5)=0$. Under $g$ we assign the weight 0 to $x_2,w_2$, 0 to $x_1,y_1,x_3, w_1,w_3,w_4$, and 1 to $x_4,x_5$ and obtain a PRDF of $T$ with at most weight $w(f)$. But $x_4$ and $x_5$ are assigned with 1, contradicting Observation \ref{PRO1}.

This completes the proof.
\end{proof}

Let $T'=T-x_2$ and let $g$ be a $\gamma_R^p$-function on $T'$.
Since $T$ is stable, $w(g)=\gamma_R^p(T-\{x_1,x_2,y_1\})+2=w(f)$.\\
If $g(x_3)=2$, then under $g$ we assign 1 to $x_1,y_1$, 0 to $x_2$ and obtain a PRDF of $T$ with same weight $w(f)$. But $x_1$ and $y_1$ are assigned with 1, contradicting Observation \ref{PRO1}.\\
If $g(x_3)=1$, then under $g$ we assign 0 to $x_1,y_1$, 2 to $x_2$ and obtain a PRDF of $T$ with same weight $w(f)$. But $x_3$ is assigned with 1, contradicting Observation \ref{PRO1}.

Thus $g(x_3)=0$ and $g(x_4)=2$.
If $f(x_5)\geq 1$, then under $g$ we assign 0 to $x_1,y_1$, 2 to $x_2$, 1 to $x_3,x_4$ and obtain a PRDF of $T$ with same weight $w(f)$. But $x_5$ is assigned with 1, contradicting Observation \ref{PRO1}.
Therefore $f(x_5)=0$. Under $g$ we assign 0 to $x_1,y_1,x_3$, 2 to $x_2$, 1 to $x_4,x_5$ and obtain a PRDF of $T$ with same weight $w(f)$. But $x_4,x_5$ is assigned with 1, contradicting Observation \ref{PRO1}.
Then the proof is complete.
\end{proof}

\begin{thm1}\label{perfiff}
{\em Let $T$ be a tree  of order $n\ge 3$. Then $T$ is a perfect Roman domination stable tree if and only if $T\in \mathcal{T}$.}
\end{thm1}
\begin{proof}
According to Lemma \ref{first}, we need only to prove necessity. Let $T$ be a perfect Roman domination stable tree of order $n\ge 3$.
The proof is by induction on $n$. If $n=3$, then $T=P_3 \in \mathcal{T}$. Let $n\ge 4$ and let the
statement hold for all perfect Roman domination stable trees of order less than $n$.
 Clearly, ${\rm diam}(T)\ge 2$.
If ${\rm diam}(T)= 2$,then $T$ must be a star of order $n\ge 4$. Clearly $T$ is not a stable tree, a contradiction.
If ${\rm diam}(T)= 3$, then $T$ is a double star. Let $T=DS_{p,q}$ be a double star with respectively $p$ and $q$ leaves attached at two support vertex. It easily proved that a stable tree is not a stable tree, a contradiction.

Consequently, we have ${\rm diam}(T)\ge 4$.
Let $P=x_1x_2\cdots x_k$ be a longest path of $T$ and $f$ be a $\gamma_R^p$-function of $T$.
According to Lemma \ref{Lemma5}, we have $d(x_2)=2$.

Now let us prove $d(x_3)=2$. Suppose to the contrary that $d(x_3)\geq 3$.
First, we prove $L(x_3)=\emptyset$. Suppose that $L(x_3)\neq\emptyset$, let $w_1 \in L(x_3)$. By Observation \ref{PRO1}, we have $f(x_1)=0$, $f(x_2)=2$, $f(w_1)\not\in \{1,2\}$.
Then we have  $f(w_1)=0$, then $f(x_3)=2$.  Now under $f$ assign  1 to  $x_1$, 0 to $x_2$ and obtain a PRDF of $T$ with fewer weight $w(f)$, a contradiction.
Therefore we have $L(x_3)=\emptyset$.

Second, we prove that $x_3$ is not adjacent to a path $P_2=w_1w_2$ with $x_3w_1 \in E(T)$ and $w_2 \in L(w_1)$. Easily we have $f(x_1)=f(w_2)=0$ and $f(x_2)=f(w_1)=2$.
It follows from Observation \ref{PRO1} and definition of $f$  that $f(x_3)=2$.
Then  under $f$ assign 1 to $x_1,w_2$, 0 to $x_2,w_1$ and obtain a PRDF of $T$ with fewer weight $w(f)$, a contradiction.


Consequently, we have $d(x_3)=2$. Hence $x_4$ must be adjacent to pendent a path $P_3=s_1s_2s_3$ for which $s_3 \in N(x_4)$.
Let $T'=T-\{s_1,s_2,s_3\}$. By Lemma \ref{AddP3}, we have $\gamma_R^p(T)=\gamma_R^p(T')+2$.
By Lemmas \ref{AddP3},\ref{AddP2} and \ref{OneV}, we have $T$ is perfect Roman domination stable iff
$T'$ is perfect Roman domination stable.

Now we show that $x_4\in W(T')$.
Let $f'$ be a $\gamma_R^p$-function of $T'$ such that $f'(x_4)\neq 0$, and we have $w(f')=\gamma_R^p(T)-2$.
If $f'(x_4)=2$, then under $f'$ we assign  1 to $s_2,s_1$ and obtain a PRDF of $T$ with same weight $w(f)$. But $s_2,s_1$ is assigned with 1, contradicting Observation \ref{PRO1}.
If $f'(x_4)=1$, then under $f'$ we assign  0 to $s_3,s_1$, 2 to $s_2$ and obtain a PRDF of $T$ with same weight $w(f)$. But $x_4$ is assigned with 1, contradicting Observation \ref{PRO1}.

Therefore, $T$ is obtained by Operation $\mathcal{O}_1$ by applying $T_i=T'$.
\end{proof}

By the construction of $\mathcal{T}$, we have
\begin{cor1}
Let $T$ be an $n$-vertex perfect Roman domination stable tree, then $n \equiv 0$ (mod 3) and $\gamma_R^p(T)=\frac{2n}{3}$.
\end{cor1}


\begin{thebibliography}{22}
\bibitem{Bauer}
D. Bauer, F. Harary, J. Nieminen, C. Suffel, Domination alternation sets in graphs, Discrete Math. 47 (1983) 153-161.

\bibitem{Chambers2009} E. W. Chambers, B. Kinnersley, N. Prince,   Extremal Problems for Roman Domination, SIAM J. Discrete Math. 23(3) (2009)  1575-1586.

\bibitem{Cockayne2004} E. J. Cockayne, P. M. Dreyer Jr., S. M. Hedetniemi, S. T. Hedetniemi, Roman domination in graphs, Discrete Math. 278 (2004)  11-22.

\bibitem{Favaron2009} O. Favaron, H. Karami, R. Khoeilar, et al, On the Roman domination number of a graph, Discrete Math. 309(10) (2009)  3447-3451.

\bibitem{Hajian}
 M. Hajian, N.J. Rad, On the Roman domination stable graphs, Discuss. Math. Graph Theory 37 (2017) 859-871.

\bibitem{Hansberg}
A. Hansberg, N. Jafari Rad and L. Volkmann, Vertex and edge critical Roman domination in graphs, Util. Math. 92 (2013) 73-88.

\bibitem{Haynes1998} T. W. Haynes, S. T. Hedetniemi, P. J. Slater (Editors), Domination in Graphs: Advanced Topics, Marcel Dekker, New York, 1998.

\bibitem{Henning}
 M.A. Henning, W.F. Klostermeyer, G. MacGillivray, Perfect Roman domination in trees, Discret. Appl. Math. (2017) In Press.

\bibitem{Rad1}
N.J. Rad, E. Sharifi, M. Krzywkowski, Domination stability in graphs, Discrete Math. 339 (7) (2016) 1909-1914.

\bibitem{Rad2}
N.J. Rad, L. Volkmann, Changing and unchanging the Roman domination number of a graph, Util. Math. 89 (2012) 79-95.

\bibitem{WW1}
S. Wang, B. Wei, A note on the independent domination number versus the domination number in bipartite graphs, Czech. Math. J., 67 (142)
(2017)  533-536.

\bibitem{WW2} M. Nandi, S. Parui, A. Adhikari, The domination numbers of cylindrical grid graphs, Appl. Math. Comput. 217 (2011) 4879-4889.
\bibitem{WW3} O. Ore, Theory of Graphs, American Mathematical Society, Providence, R.I., 1967.
\bibitem{WW4} I. Gorodezky, Domination in Kneser graphs, Dr. Sci. Thesis, University of Waterloo, Waterloo, Ontario, Canada, 2007.
\bibitem{WW5} S. Alanko, S. Crevals, A. Isopoussu, P. R. J. Ostergard, V. Pettersson, Computing the domination number of grid graphs, The Electronic Journal of Combinatorics 18 (2011) p141.
\bibitem{WW6}  P. Pavlic and J. Zerovnik, Roman domination number of the Cartesian products of paths and cycles, The Electronic Journal of Combinatorics 16 (2012) P19.
\bibitem{WW7} I. Stewart, Defend the Roman Empire!, Scientific American 281 (1999) 136-138.
\bibitem{WW8} J. Southey, M. A. Henning, Domination versus independent domination in cubic graphs, Discrete Math. 313 (2013) 1212-1220.
\bibitem{WW9} 
 J.  Liu,  X. Pan, L. Yu, D. Li, Complete characterization of bicyclic graphs with minimal Kirchhoff
 index, Discrete Appl. Math. 200 (2016) 95-107.
\bibitem{Samodivkin}
 V. Samodivkin, Roman domination in graphs: the class $\mathcal{R}_{UVR}$, Discrete Math. Algorithms Appl. 8 (2016) 1650049.
\end{thebibliography}
\end{document}